\documentclass[bj,authoryear]{imsart}

\usepackage{float,bbold,mathtools}

\startlocaldefs
\numberwithin{equation}{section}
\theoremstyle{plain}                    % default
\newtheorem{thm}{Theorem}[section]
\newtheorem{lem}{Lemma}[section]

\newtheorem{prop}{Proposition}[section]
\theoremstyle{remark}
\newtheorem{rem}{Remark}               %\renewcommand{\therem}{}

\def\N{{\mathbb N}}
\def\Z{{\mathbb Z}}
\def\R{{\mathbb R}}

\def\F{{\mathcal F}}
\def\1{{\mathbb{1}}}
\def\wtl{{\widetilde{\sigma}}}

\def\wtx{{\widetilde{X}}}
\def\wtc{{\widetilde{C}}}
\def\ninfty{\mathop{\longrightarrow}\limits_{n\to\infty}}
\def\tinfty{\mathop{\longrightarrow}\limits_{t\to\infty}}

\def\argmin{\mbox{arg}\,\min}
\newcommand{\var}{\mathop{\rm var}\nolimits}
\newcommand{\cov}{\mathop{\rm cov}\nolimits}

\newcommand{\be}{\begin{equation}}
\newcommand{\bd}{\begin{displaymath}}
\newcommand{\ed}{\end{displaymath}}
\newcommand{\bea}{\begin{eqnarray}}
\newcommand{\eea}{\end{eqnarray}}
\newcommand{\bean}{\begin{eqnarray*}}
\newcommand{\eean}{\end{eqnarray*}}

\endlocaldefs

\begin{document}

\begin{frontmatter}
\title{A log-linear model for non-stationary time series of counts}
%\title{A sample artiote that latter
\runtitle{Log-linear count processes}
%\thankstext{T1}{A sample of additional note to the title.}

\begin{aug}
\author[A]{\fnms{Anne}~\snm{Leucht}\ead[label=e1]{E-mail: anne.leucht@uni-bamberg.de}\orcid{0000-0003-3295-723X}}
\author[B]{\fnms{Michael H.}~\snm{Neumann}\ead[label=e2]{E-mail: michael.neumann@uni-jena.de}\orcid{0000-0002-5783-831X}}
%\author[B]{\inits{???}\fnms{???}~\snm{???}\ead[label=e3]{???@???}}
%%%%%%%%%%%%%%%%%%%%%%%%%%%%%%%%%%%%%%%%%%%%%%
%% Addresses                                %%
%%%%%%%%%%%%%%%%%%%%%%%%%%%%%%%%%%%%%%%%%%%%%%
\address[A]{Universit\"at Bamberg,
Institut f\"ur Statistik,
Feldkirchenstra{\ss}e 21,
D -- 96052 Bamberg,
Germany\printead[presep={,\ }]{e1}}
\address[B]{Friedrich-Schiller-Universit\"at Jena,
Institut f\"ur Mathematik,
Ernst-Abbe-Platz 2,
D -- 07743 Jena,
Germany\printead[presep={,\ }]{e2}}
%\address[B]{???\printead[presep={,\ }]{e2,e3}}
\end{aug}

\begin{abstract}
We propose a new model for non-stationary integer-valued time series which is particularly
suitable for data with a strong trend. In contrast to popular Poisson-INGARCH models,
but in line with classical GARCH models, we propose to pick the conditional distributions
from nearly scale invariant families where the mean absolute value and the standard deviation are
of the same order of magnitude.
As an important prerequisite for applications in statistics, we prove absolute regularity
of the count process with exponentially decaying coefficients.
\end{abstract}


\end{frontmatter}

%%%%%%%%%%%%%%%%%%%%%%%%%%%%%%%%%%%%%%%%%%%%%%
%%%% Main text entry area:

%%%%%%%%%%%%%%%%%%%%%%%%%%%%%%%%%%%%%%%%%%%%%%%%%%%%%%%%%%%%%%%%%%%%%%%%%%%%%%%
\section{Motivation and introduction of the model}
\label{S1}
%%%%%%%%%%%%%%%%%%%%%%%%%%%%%%%%%%%%%%%%%%%%%%%%%%%%%%%%%%%%%%%%%%%%%%%%%%%%%%%

\begin{figure}[h]
	\includegraphics[width=6.8cm]{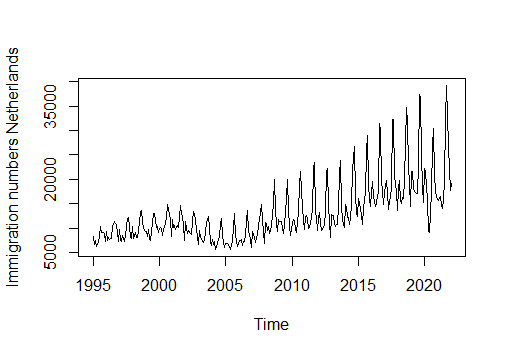}\qquad\includegraphics[width=6.8cm]{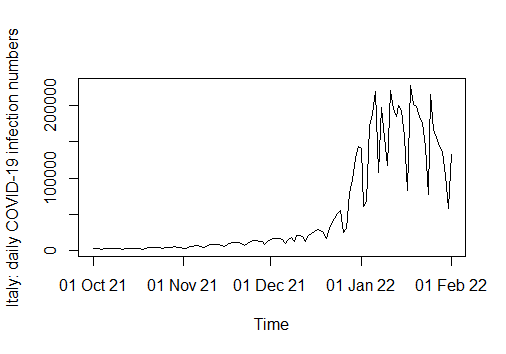} 
	\caption{left: Monthly immigration numbers for the Netherlands with increasing trend and strongly increasing seasonality;
	right: daily COVID-19 infection numbers from Italy with explosive trend.}\label{fig1}
\end{figure}

We propose a new model for time series of counts which is particularly appropriate
for modeling explosive processes. So far, the literature on models for integer-valued
time series is dominated by processes where the distribution of the count variables
conditioned on the past is taken from the family of Poisson distributions,
and where the intensities themselves are random and depend on lagged values of the count 
and the intensity variables; see e.g.~Chapter~4 in \citet{Wei18}.
While most of the results on statistical inference for these models are restricted to stationary time series,
most of real life count data exhibit strong seasonal patterns and trends, see Figure~\ref{fig1}.
These features have to be incorporated into the model as simple detrending and deseasonalization
is not feasible due to the discrete structure of the data.
With a view towards a possible explosive behavior of processes to be modeled,
we believe that it is natural that the expected values of the count variables
are of the same order of magnitude as the respective standard deviations.
This is in contrast to the family of Poisson distributions where, for $X\sim\mbox{Poi}(\sigma)$,
$EX=\sigma$ but $\sqrt{\var(X)}=\sqrt{\sigma}$, and so $\sqrt{\var(X)}/EX$ tends to zero as $\sigma\to\infty$.

Note that overdispersion in the sense of $\var(X_t)=\eta EX_t$ for some $\eta>1$ can be incorporated in
\mbox{(log-)}linear Poisson regression models by including exogenous regressors,
see e.g.~\citet{CT86} and \citet{XXGF12}. This approach is frequently used in various fields; e.g.~in demography,
health, and biology	(as in~\citet{RH03}, \citet{Petal15}, and \citet{Aagaard-etal-18})
and can be adapted to non-linear Poisson autoregressions:
If
\begin{displaymath}
X_t\mid \sigma(\sigma_0,X_0,\ldots,\sigma_{t-1},X_{t-1}) \,\sim\, \mbox{Poi}(\sigma_t\, Z_t)\qquad
\text{with }\quad \sigma_t \,=\, f(\sigma_{t-1},X_{t-1})
\end{displaymath}
where	$(Z_t)_t$ is a sequence of i.i.d.~random variables with mean~$\mu_z$ and variance~$\sigma_z^2\in(0,\infty)$
such that~$Z_t$ is stochastically independent of $\sigma(\sigma_0,X_0,\ldots,\sigma_{t-1},X_{t-1})$,
then it follows from 
\begin{eqnarray*}
E[X_t\mid \mathcal \sigma(\sigma_0,X_0,\ldots,\sigma_{t-1},X_{t-1})] & = & f(\sigma_{t-1}, X_{t-1})\, \mu_z 
		\quad\text{and}\\
\var(X_t\mid\sigma(\sigma_0,X_0,\ldots,\sigma_{t-1},X_{t-1})) & = & f^2(\sigma_{t-1}, X_{t-1})\, \sigma_z^2+f(\sigma_{t-1}, X_{t-1})\,\mu_z 
\end{eqnarray*}
that the degree of conditional dispersion can be changed adapting mean and variance of the exogenous regressors accordingly. 
Indeed, several well-known count time series models are special cases of these so-called mixed Poisson models.
For instance, if the $Z_t$'s are binomial, then the conditional distribution of~$X_t$ is a zero-inflated Poisson distribution
and $Z_t$'s being Gamma distributed result in a negative binomial distribution, see e.g.~\citet{Ketal21} and \citet{DLN22}.

Here, we propose an alternative approach to assure that conditional expectation and variance are of the same order of magnitude. We pick the conditional distributions from a family of (nearly)
scale-invariant distributions. Candidates for such distributions can be found
by discretizing scale-invariant continuous distributions. 
Let~$Y$ be a non-negative random variable with a probability density~$p$.
Then we can define  related integer-valued random variables $X_\sigma$ by setting
\begin{displaymath}
X_\sigma \,=\, \lfloor \sigma Y \rfloor,
\end{displaymath}
i.e.~$X_\sigma=k$ if and only if $\sigma Y\in[k,k+1)$.
In what follows we denote the distribution of $X_\sigma$ by~$P_\sigma$ and the corresponding
probability mass function by~$p_\sigma$. For example,
if~$Y$ is exponentially distributed with rate parameter~$1$, then 
$\sigma Y$ is exponentially distributed with rate parameter~$\sigma$ and $E[\sigma Y]=\sqrt{\var(\sigma Y)}=\sigma$.
We obtain for the corresponding integer-valued random variable~$X_\sigma$ that
\begin{displaymath}
P_\sigma\big( \{k\} \big) \,=\, e^{-k/\sigma} \,-\, e^{-(k+1)/\sigma}
\,=\, (1-p)^k\, p \qquad \forall k\in\N_0,
\end{displaymath}
where $p=1-e^{-1/\sigma}$. In this case, $X_\sigma$ has a geometric distribution with success parameter~$p$.
Furthermore, $EX_\sigma=(1-p)/p=e^{-1/\sigma}/(1-e^{-1/\sigma})$ and
$\sqrt{\var(X_\sigma)}=\sqrt{(1-p)/p^2}=e^{-1/(2\sigma)}/(1-e^{-1/\sigma})$ are of the same order of magnitude as $\sigma\to\infty$.
(This version of a geometric distribution has support~$\N_0$ and describes the number of failures
before the first success of independent Bernoulli trials with success parameter~$p$.)
Another example can be generated from the family of normal distributions.
If $Y\sim N(0,\nu^2)$, then~$|Y|$ has a so-called half-normal distribution, and corresponding integer-valued
random variables can be generated by setting $X_\sigma=\lfloor \sigma Y \rfloor$.
In both cases, the families $(P_\sigma)_{\sigma>0}$ satisfy the conditions imposed in this paper.

We impose a GARCH-type structure for the count process, i.e.
\begin{subequations}
\begin{equation}
\label{1.1a}
X_t\mid \F_{t-1} \,\stackrel{d}{=}\, \lfloor \sigma_t\, Y\rfloor \,\sim\, P_{\sigma_t},\quad t\in\N,
\end{equation}
where $\F_s=\sigma(\sigma_0,X_0,Z_0,\ldots,\sigma_s,X_s,Z_s)$ and $\sigma_t$ is a function of $\sigma_{t-1}$, $X_{t-1}$,
and an exogenous covariate $Z_{t-1}$ which may describe e.g.~seasonal effects or the effect of a changing environment.
To be specific, we will assume that $\sigma_t \,=\, f(\sigma_{t-1},X_{t-1}) \,\cdot\, Z_{t-1}$,
which can be equivalently rewritten as
\begin{equation}
\label{1.1c}
\ln(\sigma_t) \,=\, \ln( f(\sigma_{t-1},X_{t-1}) ) \,+\, C_{t-1},
\end{equation}
\end{subequations}
where $C_t=\ln(Z_t)$. Note that the initial variable $\sigma_0\geq 0$ can be chosen arbitrarily as long as some moment conditions are satisfied, see Theorem~\ref{T1} for details.

To work with such processes it is necessary to have some probabilistic properties at our disposal.
For classical GARCH processes, mixing properties have been known for a long time;
see e.g.~\citet{Bou98} for linear and \citet{CC02} and \citet{FZ06} for nonlinear variants.
These properties are typically stated for the bivariate process consisting of the observable and the
state variables.
In sharp contrast, for integer-valued GARCH, the state process $(\sigma_t)_{t\in\N_0}$ is not
mixing in general; see Remark~3 in \citet{Neu11} for a counterexample.
Moreover, mixing properties of classical GARCH processes can be deduced under weak moment conditions
on the innovation distribution (see \citet{FZ06}), while for Poisson-INGARCH time series all
conditional moments given $\sigma_t$ exist. In contrast, here we require only $E[\ln^+(Y)]<\infty$  with $\ln^+(Y)=\max\{\ln(Y), 0\}$
which means that only logarithmic conditional moments of the $X_t$ may exist. 
 
In this paper we search for conditions that allow us to prove absolute regularity ($\beta$-mixing)
of the count process $(X_t)_{t\in\N_0}\;$ ($\N_0=\N\cup\{0\}=\{0,1,2,\ldots\}$).
This will be done by a coupling approach described in greater detail in Section~\ref{S2}. 
To illustrate the usefulness of this result we discuss an application in statistics in Section~\ref{S3}.
The proof of the main result is contained in Section~\ref{S4} and the proof of asymptotic normality 
of a least squares estimator of a trend parameter is postponed to Section~\ref{S5}. Validity of a bootstrap method to quantify uncertainty of this estimator is proved in Section~\ref{S6neu}.
Proofs of a few technical results are collected in a final Section~\ref{S6}.

%%%%%%%%%%%%%%%%%%%%%%%%%%%%%%%%%%%%%%%%%%%%%%%%%%%%%%%%%%%%%%%%%%%%%%%%%%%%%%%
\section{Absolute regularity of the count process}
\label{S2}
%%%%%%%%%%%%%%%%%%%%%%%%%%%%%%%%%%%%%%%%%%%%%%%%%%%%%%%%%%%%%%%%%%%%%%%%%%%%%%%
 
Let $(\Omega,{\mathcal A},P)$ be a probability space and ${\mathcal A}_1$, ${\mathcal A}_2$
be two sub-$\sigma$-algebras of ${\mathcal A}$. Then the coefficient of absolute regularity is defined as
\bd
\beta({\mathcal A}_1,{\mathcal A}_2) \,=\, E\big[ \sup\{ |P(B\mid {\mathcal A}_1) \,-\, P(B)|\colon
\;\; B\in {\mathcal A}_2 \} \big].
\ed
For a process ${\mathbf X}=(X_t)_{t\in\N_0}$ on $(\Omega,\F,P)$, the
coefficients of absolute regularity at the point~$k$ are defined as
\bd
\beta^X(k,n) \,=\, \beta\big( \sigma(X_0,X_1,\ldots,X_k), \sigma(X_{k+n},X_{k+n+1},\ldots) \big)
\ed
and the (global) coefficients of absolute regularity as
\bd
\beta^X(n) \,=\, \sup\{ \beta^X(k,n)\colon \;\; k\in\N_0\};
\ed
see e.g.~\citet{D94}.
The intended approach to prove absolute regularity is inspired by the fact that one can construct,
on a suitable probability space $(\widetilde{\Omega},\widetilde{\F},\widetilde{P})$
two versions of the process $(X_t)_{t\in\N_0}$, $(\wtx_t)_{t\in\N_0}$ and $(\wtx_t')_{t\in\N_0}$, such that
$(\wtx_0,\ldots,\wtx_k)$ and $(\wtx_0',\ldots,\wtx_k')$ are independent and
\begin{displaymath}
\beta^X(k,n) \,=\, \widetilde{P}\left( \wtx_{k+n+r}\neq\wtx_{k+n+r}' \mbox{ for some } r\geq 0 \right).
\end{displaymath}
Since such an optimal coupling seems to be out of reach in our context we confine ourselves to
construct a ``reasonably good'' coupling.
Actually, if $(\wtx_t)_{t\in\N_0}$ and $(\wtx_t')_{t\in\N_0}$ defined on a common probability space
$(\widetilde{\Omega},\widetilde{\F},\widetilde{P})$ are any two versions of $(X_t)_{t\in\N_0}$
such that $(\wtx_0,\ldots,\wtx_k)$ and $(\wtx_0',\ldots,\wtx_k')$ are independent, then
\bea
\label{2.2}
 \beta^X(k,n)  %\nonumber \\
& \leq & \widetilde{E}\Big[ \sup_{C\in\sigma({\mathcal C})} \big\{ \big|
\widetilde{P}\left( (\wtx_{k+n},\wtx_{k+n+1},\ldots)\in C\mid \wtx_0,\ldots,\wtx_k \right) \nonumber \\
& & \qquad \qquad \qquad
\,-\, \widetilde{P}\left( (\wtx_{k+n}',\wtx_{k+n+1}',\ldots)\in C\mid \wtx_0',\ldots,\wtx_k' \right) \big| \big\} \Big]
\qquad \nonumber \\
& \leq & \widetilde{P}\left( \wtx_{k+n+r} \neq \wtx_{k+n+r}' \quad \mbox{for some } r\in\N_0 \right) \nonumber \\
& = & \widetilde{P}\left( \wtx_{k+n} \neq \wtx_{k+n}' \right) \nonumber \\
& & {} \,+\, \sum_{r=1}^\infty \widetilde{P}\left( \wtx_{k+n+r} \neq \wtx_{k+n+r}', \wtx_{k+n+r-1} = \wtx_{k+n+r-1}',
\ldots, \wtx_{k+n} = \wtx_{k+n}' \right). 
\eea
(In the second line of this display, $\sigma({\mathcal C})$ denotes the $\sigma$-algebra generated by the cylinder sets.)
We construct on a suitable probability space $(\widetilde{\Omega},\widetilde{\F},\widetilde{P})$
two versions $\big((\wtl_t,\wtx_t,\wtc_t)\big)_{t\in\N_0}$ and $\big((\wtl_t',\wtx_t',\wtc_t')\big)_{t\in\N_0}$
of the process $\big((\sigma_t,X_t,C_t)\big)_{t\in\N_0}$ such that 
\begin{equation}\label{eq.dtv}
\widetilde{P}\big( \wtx_t = \wtx_t' \mid \wtl_t,\wtl_t' \big)
\,=\, \sum_{k=0}^\infty P_{\wtl_t}(\{k\}) \wedge P_{\wtl_t'}(\{k\})
\,= \, 1 \,-\, d_{TV}\big( P_{\wtl_t}, P_{\wtl_t'} \big),
\end{equation}
where $d_{T_V}(P,Q)=(1/2)\sum_{k=0}^\infty |P({k})-Q({k})|$ denotes the total variation distance of two count measures~$P$ and~$Q$. Furthermore, if the density $p$ of $Y$ is nonincreasing on $[0,\infty)$ or differentiable everywhere on $(0,\infty)$
and $\int_0^\infty x|p'(x)|\, dx<\infty$, then Lemma~\ref{L.tv} shows that
\begin{displaymath}
d_{TV}\big( P_\sigma\,,\, P_{\sigma'} \big) \,=\, O\big( |\ln(\sigma) \,-\, \ln(\sigma')| \big).
\end{displaymath}
Hence, it is important to have the evolution of the process $\big( |\ln(\wtl_t)-\ln(\wtl_t')| \big)_{t\in\N_0}$
under control. 
Since $E\big| \ln(\lfloor \sigma Y\rfloor +1) - \ln(\lfloor \sigma' Y\rfloor +1)\big| \leq \gamma\,|\ln(\sigma)-\ln(\sigma')|$
for 
\begin{displaymath}
\gamma \,:=\, \int_0^\infty \sup\{p(y)\colon\, y\geq x\}\,dx
\end{displaymath}
(see Lemma~\ref{LA.2} below) we
suppose that the link function~$f$ in \eqref{1.1c} satisfies the following contractive condition.
\begin{subequations}
\begin{eqnarray}
\label{1.3a}
\lefteqn{ \big| \ln\big( f(\sigma, x) \big) \,-\, \ln\big( f(\sigma', x') \big) \big| } \nonumber \\
& \leq & a\, | \ln(\sigma) \,-\, \ln(\sigma') | \,+\, b\, | \ln(x+1) \,-\, \ln(x'+1) |
\qquad \forall \sigma,\sigma'>0, x,x'\in\N_0, \qquad
\end{eqnarray}
where $a,b\geq 0$ and $a+b\gamma<1$. Note that the latter condition implies validity of the classical contraction condition $a+b<1$ since it holds $\gamma\geq 1$ by definition.
Furthermore, we assume that $C_t$ is independent of $\F_{t-1}$ and $X_t$, and that
\begin{equation}
\label{1.3b}
M \,:=\, \sup_t E| C_t \,-\, E C_t | \,<\, \infty.
\end{equation}
\end{subequations}
Note that this includes the specification proposed in \citet{FT11},
\begin{displaymath}
\ln(\sigma_t) \,=\, a\, \ln(\sigma_{t-1}) \,+\, b\, \ln(X_{t-1}+1) \,+\, d \qquad \forall t\in\N,
\end{displaymath}
where the count variable~$X_t$ given the past had 
a Poisson distribution with random intensity~$\sigma_t$, and $d$ was a constant.
In that case, for appropriately chosen values of the parameters~$a$ and~$b$, there exists a stationary
process with such a dynamics.
Lemmas~\ref{L4.1},~\ref{L4.2}, and~\ref{L.tv}
allow us to prove the following main result of this contribution.

{\thm
\label{T1}
Suppose that \eqref{1.1a}, \eqref{1.1c}, \eqref{1.3a}, and \eqref{1.3b} are fulfilled, and let
$E|\ln(\sigma_0)|+E\ln^+(Y)<\infty$.
We assume that the density $p\colon [0,\infty)\rightarrow [0,\infty)$ of~$Y$ is
\begin{itemize}
\item[(i)] monotonously non-increasing\\
or
\item[(ii)] everywhere differentiable.
\end{itemize}
In the former case we set $\gamma=\Gamma=1$ and in the latter $\gamma=\int_0^\infty \sup\{p(y)\colon\, y\geq x\}\, dx$
and $\Gamma=\big( 1+\int_0^\infty x\,|p'(x)|\,dx \big)/2$.
Then the count process $(X_t)_{t\in\N_0}$ is absolutely regular ($\beta$-mixing), and the corresponding coefficients satisfy
\begin{displaymath}
\beta^X(n) \,\leq\, (a+b\gamma)^n\, \frac{\Gamma}{1-a} \, \Big\{ 2\, E|\ln(\sigma_0)|
\,+\, \frac{2b\,(\|p\|_\infty + E\ln^+(Y)) \,+\, 2M}{1-a-b} \Big\}.
\end{displaymath}
}

\noindent
Theorem~\ref{T1} can serve as a basis for various statistical applications.
For instance, confidence sets and statistical tests can be developed relying on \citeauthor{R95}'s \citeyearpar{R95}
CLT for triangular arrays of nonstationary random variables; see Section~\ref{S3} for details.
 
{\rem The monotonicity assumption on~$p$ is satisfied for instance if~$Y$ is exponentially distributed or if it is half-normal.
If~$Y$ has a chi-square distribution with at least $k=3$ degrees of freedom, then the corresponding
density~$p_k$ is unimodal with a mode at~$k-2$. However it is differentiable everywhere on $(0,\infty)$ and
\begin{eqnarray*}
\int_0^\infty x\, |p_k'(x)| \, dx
& = & \int_0^{k-2} x\, p_k'(x) \, dx \;-\; \int_{k-2}^\infty x\, p_k'(x) \, dx \\
& = & 2(k-2)p_k(k-2) \,-\, \int_0^{k-2} p_k(x)\, dx \,+\, \int_{k-2}^\infty p_k(x) \, dx \,<\, \infty.
\end{eqnarray*}
Hence, our conditions on the distribution of~$Y$ are satisfied.
Another example is given by the family of Cauchy distributions.
Such a distribution has a density $p(x)=\pi \gamma/((x-\mu)^2+\gamma^2)$
on~$\R$, where $\mu\in\R$ is the location and $\gamma>0$ the scale parameter.
If~$Z$ follows such a distribution we would call the distribution of $Y=|Z|$ to be half-Cauchy.
Such a distribution does not have finite moments of order greater than or equal to one, however, $E\ln^+(Y)$ is finite.
$Y$ has a density $p(x)=\pi \gamma/((x-\mu)^2+\gamma^2)+\pi \gamma/((x+\mu)^2+\gamma^2)$ which is not non-increasing
on $[0,\infty)$ if $|\mu|$ is large, however, it is decreasing on $[|\mu|,\infty)$.
Again, our conditions are satisfied by this family of distributions.
}

{\rem The proof of Theorem~\ref{T1} relies heavily on the discretization $\lfloor \sigma Y\rfloor$ to bound \eqref{2.2} and \eqref{eq.dtv},
respectively. Still, these two relations build the basis for proving absolute regularity of other count time series such as
Poisson INGARCH processes and variants thereof, see e.g.~\citet{DN19, Neu21}, or \citet{DLN22}.}

%%%%%%%%%%%%%%%%%%%%%%%%%%%%%%%%%%%%%%%%%%%%%%%%%%%%%%%%%%%%%%%%%%%%%%%%%%%%%%%
\section{An application in statistics}
\label{S3}
%%%%%%%%%%%%%%%%%%%%%%%%%%%%%%%%%%%%%%%%%%%%%%%%%%%%%%%%%%%%%%%%%%%%%%%%%%%%%%%

Suppose that we observe the random variables $X_1,\ldots,X_n$ which follow a log-linear
model similar to that proposed in \citet{FT11}.
However, we include an increasing intercept term which causes a strong trend.
To be specific, we assume that
\begin{equation}
\label{eq3.0}
\ln(\sigma_t) \,=\, a\, \ln(\sigma_{t-1}) \,+\, b\, \ln(X_{t-1}+1) \,+\, c\, \ln(t),
\end{equation}
where $a\geq 0$, $b,c>0$, and $a+b\gamma<1$. We also assume that
\begin{displaymath}
	X_t \mid \F_{t-1} \,\stackrel{d}{=}\, \lfloor \sigma_t Y \rfloor.
\end{displaymath}
Although consistent estimation of the parameters $a$, $b$ and $c$ is a highly relevant issue, it is left for future research as it is mathematically demanding. On the one hand, the $\sigma_t$'s cannot be observed. Hence,  one  has to come up with  suitable estimators for the volatility process as a preliminary step if $a\neq 0$. Moreover, we expect the convergence rate of least-squares parameter estimators to be comparatively poor in view of the multicollinarity of
$(\ln(t))_{t=2,\dots, n}$ and $(\ln(X_{t-1}+1))_{t=2,\dots, n}$.
However, under mild regularity conditions we obtain that
\begin{equation}
\label{eq3.1}
E \ln(X_t+1) \,=\, \frac{c}{1-a-b} \ln(t) \,+\, O(1);
\end{equation}
see (\ref{eq5.1}) and (\ref{eq5.2}) in Section~5 below.
Hence, the parameter $\theta:=c/(1-a-b)$ characterizes the trend in the data.
To estimate it, we may fit the regression model
\begin{displaymath}
\ln( X_t+1 ) \,=\, \theta\, \ln(t) \,+\, \varepsilon_t, \qquad t=1,\ldots,n.
\end{displaymath}
The ordinary least squares estimator is given as
\begin{displaymath}
\widehat{\theta}_n  \,=\, \frac{ \sum_{t=1}^n \ln(t) \, \ln(X_t+1) }{ \sum_{t=1}^n \big( \ln(t) \big)^2 }
\end{displaymath}
and it is centered about the best projection $\bar{\theta}=\bar{\theta}(n)$ which is given as
\begin{displaymath}
\bar{\theta}
\,=\, \argmin_\theta \sum_{t=1}^n E\big| \ln(X_t+1) \,-\, \theta\, \ln(t) \big|^2
\,=\, \frac{ \sum_{t=1}^n \ln(t) \, E\ln(X_t+1) }{ \sum_{t=1}^n \big( \ln(t) \big)^2 }.
\end{displaymath}

Using the mixing property derived in Theorem~\ref{T1} in conjunction with \citeauthor{R95}'s \citeyearpar{R95} 
central limit theorem for nonstationary and mixing random variables we can prove the following central limit theorem.
In particular, it shows a non-standard rate of convergence of our estimator.

{\prop
\label{P3.1}
Suppose that the conditions of Theorem~\ref{T1} are satisfied and
\begin{displaymath}
\ln(\sigma_t) \,=\, a\, \ln(\sigma_{t-1}) \,+\, b\, \ln(X_{t-1}+1) \,+\, c\, \ln(t),
\end{displaymath}
where $a\geq 0$, $b,c>0$, and $a+b\gamma<1$. Moreover, assume that $\sigma_0\geq 1$ and
$\|\ln(\sigma_0)\|_{2+\delta}+\|\ln^+(Y)\|_{2+\delta}<\infty$ for some $\delta>0$. Then, 
\begin{displaymath}
\sqrt{n}\, \ln(n)\, \big( \widehat{\theta}_n \,-\, \bar{\theta} \big)
\,\stackrel{d}{\longrightarrow}\, \mathcal N\big( 0,\, \sigma^2\big),
\end{displaymath}
where $\sigma^2=\var\big( \ln(Y) \big) \, (1-a)^2/(1-(a+b))^2$.
}

{\rem
\label{R3}
It is known that $\int_0^1 |\ln(y)|^k\, dy=k!$ $\forall k\in\N$.
Therefore, $\int_0^1 |\ln(y)|^{2+\delta} p(y)\, dy<\infty$ if $\|p\|_\infty<\infty$, which implies that
$E\big[|\ln(Y)|^{2+\delta}\big]\leq \|p\|_\infty \, \int_0^1 |\ln(y)|^{2+\delta} \, dy+E\big[|\ln^+(Y)|^{2+\delta}\big]<\infty$.
This means in particular that the variance of~$\ln(Y)$ is finite under the above conditions.
}

\begin{figure}[h]
	\includegraphics[width=7cm]{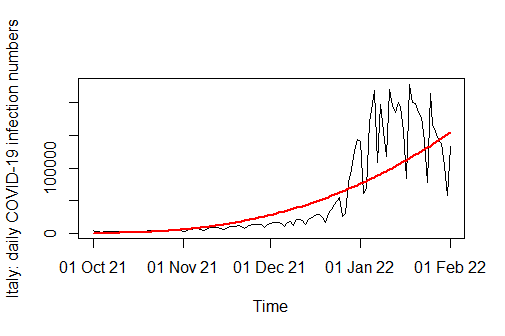} 
	\caption{COVID-19 infection numbers with estimated trend curve.}\label{fig3}
\end{figure}

\begin{figure}[h]
\includegraphics[width=0.9\textwidth]{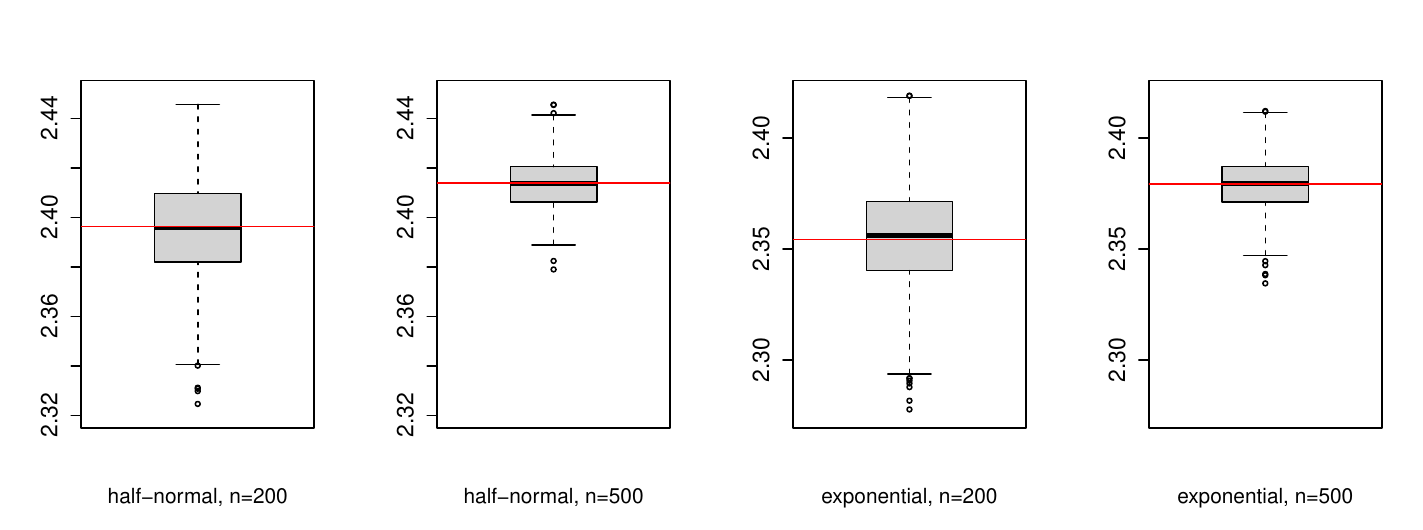}
\caption{Boxplots for sample sizes $n=200$ and $n=500$: left: $Y\sim$ half-normal  with $EY=1$, right: $Y\sim Exp(1)$. }\label{fig2}
\end{figure}

To illustrate the dynamics of these models, we revisit the COVID-19 data from Italy. We estimate $\widehat\theta_n=2.48$
and visualize the curve $t\mapsto t^{\widehat\theta_n}$, see Figure~\ref{fig3}.
Moreover, we investigate the finite sample behaviour of our estimator in synthetic samples of sizes $n=200$ and $n=500$
with parameters chosen as $a=b=0.1$ and $c=2$ which give $\theta=2.5$,
see Figure~\ref{fig2}. The distribution of~$Y$ is chosen to be either half-normal (left) or exponential (right), both with expectation~1. The box plots are based on 1000 iterations and the red line marks the target value $\bar\theta$. We obtained
$\bar\theta$ by simulation with 20000 Monte Carlo loops as it depends on the unknown distribution of the $X_t$'s as well as the sample size. The performance of the proposed estimator is very convincing in the sense that the median is very close to the target value also for smaller sample size and that the variability becomes smaller as sample sizes increase.

It is possible to construct an asymptotic confidence interval for the parameter~$\bar{\theta}$.
To this end, we have to approximate the distribution of
\begin{displaymath}
T_n \,:=\, \sqrt{n}\, \ln(n)\, \big(\widehat{\theta}_n-\bar{\theta}\big)
\,=\, \sum_{t=1}^n w_t\, \big[ \ln(X_t+1) \,-\, E\ln(X_t+1) \big],
\end{displaymath}
where $w_t=\sqrt{n}\, \ln(n)\, \ln(t)/\sum_{s=1}^n (\ln(s))^2$.
This can be conveniently achieved using the so-called dependent wild bootstrap
which was introduced by \citet{Sha10} for smooth functions of the mean.
In our context, we approximate the distribution of~$T_n$ by
the conditional distribution of
\begin{displaymath}
T_n^* \,:=\, \sum_{t=1}^n w_t\, \big[ \ln(X_t+1) \,-\, \widehat{m}_n(t) \big]\, W_t^*,
\end{displaymath}
where $\widehat{m}_n(t)$ is an appropriate estimator of $m(t)=E\ln(X_t+1)$,
and $(W_t^*)_{t=1,\ldots,n}$ is a triangular array of auxiliary random variables $W_t^*=W_{t,n}^*$
which are independent of the original sample $X_1,\ldots,X_n$, satisfy
$EW_t^*=0$ and $\cov(W_s^*,W_t^*)=\rho(|s-t|/l_n)$, where $\rho(u)\rightarrow_{u\to 0}1$
and $\sum_{r=1}^{n-1}|\rho(|r|/l_n)|=O(l_n)$.
The parameter~$l_n$ plays a similar role as the block length in blockwise bootstrap methods
and it has to be chosen such that $l_n\rightarrow_{n\to\infty}\infty$ and $l_n=o(n)$.
To simplify matters we suppose that $W_1^*,\ldots,W_n^*$ are jointly Gaussian.
It turns out that the function~$m$ is monotonously increasing (see Lemma~\ref{L.monotone} below) and that
$m(t)=O(\ln(t))$. Hence $m(t)$ can be well estimated by a nearest neighbor estimator
\begin{displaymath}
\widehat{m}_n(t) \,=\, \sum_{s\colon |s-t|\leq N_n} \ln(X_s+1)/\#\{s\in\{1,\ldots,n\}\colon\, |s-t|\leq N_n\},
\end{displaymath}
where $N_n\rightarrow\infty$ as $n\to\infty$.
It follows from (\ref{eq5.4}) and (\ref{eq5.5}) that
\begin{equation}
\label{eq3.3}
\var\big( \widehat{m}_n(t) \big) \,=\, O\big( 1/N_n \big).
\end{equation}
Furthermore, since 
$\big| E\widehat{m}_n(t) \,-\, m(t) \big| \,\leq\, E\ln(X_{(t+N_n)\wedge n}+1) \,-\, E\ln(X_{(t-N_n)\vee 1}+1)$
we obtain that
\begin{equation}
\label{eq3.4}
\sum_{t=1}^n \big| E\widehat{m}_n(t) \,-\, m(t) \big|^2 \,=\, O\big( N_n\, (\ln(n))^2 \big).
\end{equation}

{\rem
A simple way to construct the random variables $W_1^*,\ldots,W_n^*$ is to take first an Ornstein-Uhlenbeck process
$(U_t)_{t\geq 0}$, i.e.~a Gaussian process with continuous sample paths, $EU_t=0$ and $\cov(U_s,U_t)=\exp(-|s-t|)$,
and to define $W_t^*=U_{t/l_n}$, $t=1,\ldots,n$.
Then the practical implementation becomes easy since a discrete sample of an Ornstein-Uhlenbeck process 
forms an AR$(1)$ process, i.e.
\begin{displaymath}
W_t^* \,=\, e^{-1/l_n}\, W_{t-1}^* \,+\, \sqrt{1 \,-\, e^{-2/l_n}} \, \varepsilon_t^*,
\end{displaymath}
where $W_1^*,\varepsilon_2^*,\ldots,\varepsilon_n^*$ are independent standard normal variables. 
}

{\prop
\label{P3.2}
Suppose that the conditions of Proposition~\ref{P3.1} are satisfied. Furthermore, we assume that $\sigma_0=1$,
$\|\ln^+(Y)\|_{4+\delta}<\infty$ for some $\delta>0$, and that
$l_n\ninfty\infty$, $l_n/n\ninfty 0$, and $l_n/N_n+l_n N_n ((\ln(n))^2/n\ninfty 0$. Then
\begin{itemize}
\item[(i)\quad]
$\sup_{x\in\R} \big| P\big( T_n\leq x \big) \,-\, P\big( T_n^*\leq x \mid X_1,\ldots,X_n \big) \big|
\,\stackrel{P}{\longrightarrow}\, 0$.
\item[(ii)\quad]
Let $u^*_{1-\alpha/2}$ be the $(1-\alpha/2)$-quantile of $P^{T_n^*\mid X_1,\ldots,X_n}$. Then
\begin{displaymath}
P\Big( \bar{\theta} \in \big[ \widehat{\theta}_n \,-\, \frac{u^*_{1-\alpha/2}}{\sqrt{n}\, \ln(n)},
\widehat{\theta}_n \,+\, \frac{u^*_{1-\alpha/2}}{\sqrt{n}\, \ln(n)} \big] \Big) \,\ninfty\, 1 \,-\, \alpha.
\end{displaymath}
\end{itemize}
}

We illustrate the finite sample performance of the proposed bootstrap method by some simulations regarding the coverage of the bootstrap-based $(1-\alpha)$-confidence intervals for $\alpha=0.05$ and $\alpha=0.1$.  The data are generated in the same way as in the previous simulations. We use 1000 Monte Carlo loops, each consisting of 1000 bootstrap replications while varying the window size $N_n$ of the nearest-neighbor estimator and the tuning parameter $l_n$ of the dependent wild bootstrap, see Table~\ref{table.1}. We observe that the performance of our bootstrap method is convincingly stable with respect to the choice of $l_n$ as long as $N_n$ is not chosen too small. Conversely, choosing $N_n$ too large results in a higher variance of the bootstrap estimator which then gives wider confidence intervals with high coverage. 
\begin{table}\label{table.1}\caption{Rejection frequencies}
\begin{center}
\begin{tabular}{|c|c||c|c||c|c|}
	\hline
	\multicolumn{2}{|c||}{ $n=500$ }& \multicolumn{2}{|c||}{ half-normal} & \multicolumn{2}{|c|}{exponential}    \\ 
	$l_n$& $N_n$& $\alpha=0.1$ & $\alpha=0.05$   &$\alpha=0.1$ & $\alpha=0.05$       \\ \hline \hline
	20& 55& 0.900 &  0.953 & 0.860  &  0.919 \\ \hline
	20& 65 &  0.940 & 0.969 & 0.904  &  0.952 \\ \hline
	25 & 60& 0.912 & 0.953  & 0.887  &   0.933\\ \hline
	25 & 70& 0.944 &  0.976 &  0.920  & 0.956  \\ \hline
	30& 65 & 0.922  & 0.959 & 0.895  & 0.940   \\ \hline
	30& 75& 0.951 & 0.982 &  0.926 & 0.961  \\ \hline
\end{tabular}
\end{center}
\end{table}

%%%%%%%%%%%%%%%%%%%%%%%%%%%%%%%%%%%%%%%%%%%%%%%%%%%%%%%%%%%%%%%%%%%%|%%%%%%%%%%%
\section{Proof of Theorem~\ref{T1}}
\label{S4}
%%%%%%%%%%%%%%%%%%%%%%%%%%%%%%%%%%%%%%%%%%%%%%%%%%%%%%%%%%%%%%%%%%%%%%%%%%%%%%%

The proof of our main result is based on the following two lemmas.

{\lem
\label{L4.1}
Suppose that \eqref{1.1a}, \eqref{1.1c}, \eqref{1.3a}, and \eqref{1.3b} are fulfilled,
that the density~$p$ of $Y$ is bounded, and let $E|\ln(\sigma_0)|+E\ln^+(Y)<\infty$.
Let $\big((\wtl_t,\wtx_t,\wtc_t)\big)_{t\in\N_0}$ and $\big((\wtl_t',\wtx_t',\wtc_t')\big)_{t\in\N_0}$
be independent versions of $\big((\sigma_t,X_t,C_t)\big)_{t\in\N_0}$ which are defined on 
a suitable probability space $\big(\widetilde{\Omega},\widetilde{\F},\widetilde{P}\big)$. Then
\begin{displaymath}
\sup_{k\in\N_0} \widetilde{E}\big| \ln(\wtl_k) \,-\, \ln(\wtl_k') \big|
\,\leq\, 2\, E|\ln(\sigma_0)| \,+\, \frac{{2b\, (\|p\|_\infty+ E\ln^+(Y)) \,+\, 2M}}{1-a-b}.
\end{displaymath}
}

\begin{proof}
It follows from \eqref{1.1c} and \eqref{1.3a} that
\begin{eqnarray*}
\lefteqn{ \big| \ln(\wtl_t) \,-\, \ln(\wtl_t') \big| } \\
& \leq & a\, \big| \ln(\wtl_{t-1}) \,-\, \ln(\wtl_{t-1}') \big|
\,+\, b\, \big| \ln(\wtx_{t-1}+1) \,-\, \ln(\wtx_{t-1}'+1) \big|
\,+\, \big| \wtc_{t-1} \,-\, \wtc_{t-1}' \big|.
\end{eqnarray*}
Note we obtain from Lemma~\ref{LA.1} that
$\widetilde{E}\big(|\ln(\wtx_{t-1}+1)-\ln(\wtl_{t-1}+1)|\big|\wtl_{t-1},\wtl_{t-1}'\big)
=\widetilde{E}\big(|\ln(\wtx_{t-1}+1)-\ln(\wtl_{t-1}+1)|\big|\wtl_{t-1}\big)\leq  \|p\|_\infty+E\ln^+(Y)$
and, analogously,
$\widetilde{E}\big(|\ln(\wtx_{t-1}'+1)-\ln(\wtl_{t-1}'+1)|\big|\wtl_{t-1},\wtl_{t-1}'\big)\leq  \|p\|_\infty+E\ln^+(Y)$.
Therefore we obtain, for $t\geq 1$,
\begin{eqnarray*}
\lefteqn{ \widetilde{E}\Big( \big| \ln(\wtl_t) \,-\, \ln(\wtl_t') \big| \;\Big|\; \wtl_{t-1},\wtl_{t-1}' \Big) } \\
& \leq & a\, \big| \ln(\wtl_{t-1}) \,-\, \ln(\wtl_{t-1}') \big|
\,+\, b\, \big| \ln(\wtl_{t-1}+1) \,-\, \ln(\wtl_{t-1}'+1) \big| \\
& & \,+\, b\, \widetilde{E}\Big( \big| \ln(\wtx_{t-1}+1)-\ln(\wtl_{t-1}+1) \big| \;\Big|\; \wtl_{t-1},\wtl_{t-1}' \Big) \\
& & \,+\, b\, \widetilde{E}\Big( \big| \ln(\wtx_{t-1}'+1)-\ln(\wtl_{t-1}'+1) \big| \;\Big|\; \wtl_{t-1},\wtl_{t-1}' \Big) \\
& & \,+\, \widetilde{E}\Big( \big| \wtc_{t-1} \,-\, EC_{t-1} \big| \;\Big|\; \wtl_{t-1},\wtl_{t-1}' \Big)
\,+\, \widetilde{E}\Big( \big| \wtc_{t-1}' \,-\, EC_{t-1} \big| \;\Big|\; \wtl_{t-1},\wtl_{t-1}' \Big) \\
& \leq & (a+b)\, \big| \ln(\wtl_{t-1}) \,-\, \ln(\wtl_{t-1}') \big| \,+\, \bar{M},
\end{eqnarray*}
where $\bar{M}=2b\,(\|p\|_\infty+E\ln^+(Y))+2M$.
Taking expectation on both sides of this inequality we obtain that
\begin{eqnarray*}
\lefteqn{ \widetilde{E} \big| \ln(\wtl_t) \,-\, \ln(\wtl_t') \big| } \\
& \leq & (a+b) \, \widetilde{E} \big| \ln(\wtl_{t-1}) \,-\, \ln(\wtl_{t-1}') \big| \;+\; \bar{M} \\
& \leq & (a+b) \, \Big\{ (a+b)\, \widetilde{E} \big| \ln(\wtl_{t-2}) \,-\, \ln(\wtl_{t-2}') \big| \;+\; 
\bar{M} \Big\}
\;+\; \bar{M} \\
& \leq & \,\ldots\,\leq\, (a+b)^t\; \widetilde{E} \big| \ln(\wtl_0) \,-\, \ln(\wtl_0') \big|
\;+\; \bar{M} \, \big\{ 1 \,+\, (a+b) \,+\, \cdots \,+\, (a+b)^{t-1} \big\} \\
\lefteqn{ \;\leq\,  2\, E|\ln(\sigma_0)| \,+\, \frac{{\bar{M}}}{1-a-b}. \hspace*{8cm}\qedhere}
\end{eqnarray*}
\end{proof}

{\lem
\label{L4.2}
Suppose that \eqref{1.1a},  \eqref{1.1c}, \eqref{1.3a}, and \eqref{1.3b} are fulfilled, and let $E|\ln(\sigma_0)|<\infty$.
Furthermore, we assume that $p\colon [0,\infty)\rightarrow [0,\infty)$ is continuous and
$\gamma:=\int_0^\infty \sup\{p(y)\colon\, y\geq x\}\, dx<\infty$.
\\
Then there exist versions $\big((\wtl_t,\wtx_t,\wtc_t)\big)_{t\in\N_0}$ and $\big((\wtl_t',\wtx_t',\wtc_{t}')\big)_{t\in\N_0}$
of the process $\big((\sigma_t,X_t,C_t)\big)_{t\in\N_0}$ such that, for all $k,n\geq 0$,
\begin{eqnarray*}
& \mbox{\rm(i)} & \quad
\widetilde{E} \big| \ln(\wtl_{k+n}) \,-\, \ln(\wtl_{k+n}') \big|
\,\leq\, (a+b\gamma)^n \, \widetilde{E} \big| \ln(\wtl_k) \,-\, \ln(\wtl_k') \big|, \\
& \mbox{\rm(ii)} & \quad 
\widetilde{E} \Big[ \big| \ln(\wtl_{k+n+r}) \,-\, \ln(\wtl_{k+n+r}') \big|
\, \1\big( \wtx_{k+n}=\wtx_{k+n}',\ldots,\wtx_{k+n+r-1}=\wtx_{k+n+r-1}' \big) \Big] \\
& & \qquad \qquad \qquad \qquad \qquad \qquad \qquad
\leq\, a^r \, (a+b\gamma)^n \, \widetilde{E} \big| \ln(\wtl_k) \,-\, \ln(\wtl_k') \big| \qquad \forall r\geq 1, \\
& \mbox{\rm(iii)} & \quad
\widetilde{P} \big( \wtx_{k+n+r}\neq\wtx_{k+n+r}',\wtx_{k+n+r-1}=\wtx_{k+n+r-1}',\ldots,\wtx_{k+n}=\wtx_{k+n}' \big) \\
& & \qquad \leq\, \Gamma\,\widetilde{E} \Big[ \big| \ln(\wtl_{k+n+r}) \,-\, \ln(\wtl_{k+n+r}') \big| 
\, \1\big( \wtx_{k+n}=\wtx_{k+n}',\ldots,\wtx_{k+n+r-1}=\wtx_{k+n+r-1}' \big) \Big]\;\;\forall r\geq 0,
\end{eqnarray*}
where $\Gamma=1$ if the density $p$ is non-increasing on $[0,\infty)$
and $\Gamma=\big(1+\int_0^\infty x|p'(x)|\,dx\big)/2$ if $p$ is differentiable everywhere on $(0,\infty)$
with $\int_0^\infty x|p'(x)|\,dx<\infty$.
}

\begin{proof}
Let $t\geq k$. For given $\wtl_t$ and $\wtl_t'$, we apply a {\em maximal coupling} of the respective random variables
$\wtx_t$ and $\wtx_t'$, i.e.~$\wtx_t$ and $\wtx_t'$ are defined such that
\begin{itemize}
\item[a)] $\widetilde{P}\big(\wtx_t=k\mid \wtl_t,\wtl_t'\big)\,=\,P_{\wtl_t}(\{k\})
\quad\mbox{and}\quad \widetilde{P}\big(\wtx_t'=k\mid \wtl_t,\wtl_t'\big)\,=\,P_{\wtl_t'}(\{k\}) \qquad\forall k\in\N_0$,
\item[b)] $\widetilde{P}\big(\wtx_t=\wtx_t'\mid \wtl_t,\wtl_t'\big)
\,=\, 1 \,-\, d_{TV}\big(P_{\wtl_t},P_{\wtl_t'}\big)
\,=\, \sum_{k=0}^\infty P_{\wtl_t}(\{k\})\wedge P_{\wtl_t'}(\{k\})$.
\end{itemize}
(Note that our definition of the total variation norm differs from that in \citet{Lin92} by the factor~2.)
Furthermore, and in contrast to the construction used for the proof of Theorem~5.2 in \citet[Chapter~I]{Lin92},
we couple $\wtx_t$ and $\wtx_t'$ such that $\wtx_t\geq\wtx_t'$ holds with conditional probability~1
if $\wtl_t\geq\wtl_t'$ and, vice versa, $\wtx_t\leq\wtx_t'$ holds with conditional probability~1
if $\wtl_t\leq\wtl_t'$. (This is possible since $\lfloor\sigma Y\rfloor$ is stochastically greater than
$\lfloor\sigma' Y\rfloor$ if $\sigma\geq\sigma'$ and
$P^{\lfloor\sigma Y\rfloor}((-\infty,t])\leq P^{\lfloor\sigma' Y\rfloor}((-\infty,t])$ $\forall t$ implies that
$(P^{\lfloor\sigma Y\rfloor}-P^{\lfloor\sigma Y\rfloor}\wedge P^{\lfloor\sigma' Y\rfloor})((-\infty,t])
\leq (P^{\lfloor\sigma' Y\rfloor}-P^{\lfloor\sigma Y\rfloor}\wedge P^{\lfloor\sigma' Y\rfloor})((-\infty,t])$ $\forall t$.)
And finally, we choose the exogenous variables such that $\wtc_t=\wtc_t'$.

\mbox{(i)\quad}
Let $k<t\leq k+n$.
Since $\wtx_t-\wtx_t'$ has the same sign as $\wtl_t-\wtl_t'$,
$\ln(\wtx_t+1)-\ln(\wtx_t'+1)$ is either non-negative or non-positive with conditional probability~1,
and we obtain by Lemma~\ref{LA.2}
\begin{eqnarray*}
\widetilde{E} \Big( \big| \ln(\wtx_t+1) \,-\, \ln(\wtx_t'+1) \big| \;\Big|\; \wtl_{t}, \wtl_{t}' \Big)
& \leq & \gamma \, \big| \ln(\wtl_t) \,-\, \ln(\wtl_t') \big|.
\end{eqnarray*}
This yields that
\begin{displaymath}
\widetilde{E} \Big( \big| \ln(\wtl_t) \,-\, \ln(\wtl_t') \big| \;\Big|\; \wtl_{t-1}, \wtl_{t-1}' \Big) 
\,\leq \, (a+b\gamma)\, \big| \ln(\wtl_{t-1}) \,-\, \ln(\wtl_{t-1}') \big|.
\end{displaymath}
Taking expectation on both sides of this inequality, and using the resulting inequality~$n$ times we obtain
\begin{displaymath}
\widetilde{E} \big| \ln(\wtl_{k+n}) \,-\, \ln(\wtl_{k+n}') \big|
\,\leq\, (a+b\gamma)^n \, \widetilde{E} \big| \ln(\wtl_k) \,-\, \ln(\wtl_k') \big|.
\end{displaymath}
\mbox{(ii)\quad}
Let now $t=k+n+r>k+n$. Then we obtain immediately
\begin{eqnarray*}
\lefteqn{ \big| \ln(\wtl_t) \,-\, \ln(\wtl_t') \big|\, \1\big( \wtx_{k+n}=\wtx_{k+n}',\ldots,\wtx_{t-1}=\wtx_{t-1}' \big) } \\
& \leq & a\, \big| \ln(\wtl_{t-1}) \,-\, \ln(\wtl_{t-1}') \big|\, \1\big( \wtx_{k+n}=\wtx_{k+n}',\ldots,\wtx_{t-2}=\wtx_{t-2}' \big) \\
& \leq & \cdots \,\leq\, a^r\, \big| \ln(\wtl_{k+n}) \,-\, \ln(\wtl_{k+n}') \big|,
\end{eqnarray*}
which yields (ii). Statement (iii) follows from Lemma~\ref{L.tv} since b) gives
\begin{eqnarray*}
& \widetilde{P} & \!\big( \wtx_{k+n+r}\neq\wtx_{k+n+r}',\wtx_{k+n+r-1}=\wtx_{k+n+r-1}',\ldots,\wtx_{k+n}=\wtx_{k+n}' \big) 
\hspace*{1.5cm}\\
\lefteqn{\;=\, E\big(d_{TV}(\widetilde P_{\wtl_{k+n+r}}, \widetilde P_{\wtl_{k+n+r}'}) \,
\1\big( \wtx_{k+n}=\wtx_{k+n}',\ldots,\wtx_{k+n+r-1}=\wtx_{k+n+r-1}' \big)\big). \hspace*{1.3cm}\qedhere}
\end{eqnarray*}
\end{proof}

\noindent
Now we are in a position to prove our main result.

\begin{proof}[Proof of Theorem~\ref{T1}]
Let $\big((\wtl_t,\wtx_t,\wtc_t)\big)_{t\in\N_0}$ and $\big((\wtl_t',\wtx_t',\wtc_t')\big)_{t\in\N_0}$
be two versions of the process $\big((\sigma_t,X_t,C_t)\big)_{t\in\N_0}$, where 
$(\wtl_0,\wtx_0,\wtc_0,\ldots,\wtl_{k-1},\wtx_{k-1},\wtc_{k-1},\wtl_k,\wtx_k)$ and\\
$(\wtl_0',\wtx_0',\wtc_0',\ldots,\wtl_{k-1}',\wtx_{k-1}',\wtc_{k-1}',\wtl_k',\wtx_k')$
are independent, and where $\wtc_k,\wtl_{k+1},\wtx_{k+1},\wtl_{k+2},\wtx_{k+2},\ldots$
are coupled with their respective counterparts $\wtc_k',\wtl_{k+1}',\wtx_{k+1}',\wtl_{k+2}',\wtx_{k+2}',\ldots$
as described in the proof of Lemma~\ref{L4.2}.
Then we obtain from \eqref{2.2} and Lemmas~\ref{L4.1} and~\ref{L4.2} that
\begin{eqnarray*}
\lefteqn{ \beta^X(k,n) } \nonumber \\
& \leq & \widetilde{P}\left( \wtx_{k+n} \neq \wtx_{k+n}' \right) \\
& & {} \,+\, \sum_{r=1}^\infty \widetilde{P}\left( \wtx_{k+n+r} \neq \wtx_{k+n+r}', \wtx_{k+n+r-1} = \wtx_{k+n+r-1}',
\ldots, \wtx_{k+n} = \wtx_{k+n}' \right) \\
\lefteqn{ \;\leq\, \Gamma\, \sum_{r=0}^\infty a^r \, (a+b\gamma)^n\, \big\{ 2\, E\big|\ln(\sigma_0)\big|
\,+\, \frac{{2b\,(\|p\|_\infty+E\ln^+(Y)) \,+\, 2M}}{1-a-b} \big\}. \hspace*{1.85cm}\qedhere}
\end{eqnarray*}
\end{proof}

%%%%%%%%%%%%%%%%%%%%%%%%%%%%%%%%%%%%%%%%%%%%%%%%%%%%%%%%%%%%%%%%%%%%%%%%%%%%%%%
\section{Proof of Proposition~\ref{P3.1}}
\label{S5}
%%%%%%%%%%%%%%%%%%%%%%%%%%%%%%%%%%%%%%%%%%%%%%%%%%%%%%%%%%%%%%%%%%%%%%%%%%%%%%%

The following lemma is a consequence of a CLT in \citet{R95} and it serves as a basis for the proof
of Proposition~\ref{P3.1}.

{\lem
\label{p.clt}
Suppose that $(X_{t,n})_{t=1,\dots n},\; n\in\N,$ is a triangular array of centered random variables with
$\sup_{t\leq n}(\|X_{t,n}\|_{2+\delta})=O(n^{-1/2})$ 
for some $\delta>0$, where  $\|X_{t,n}\|_q:=(E|X_{t,n}|^q)^{1/q}$. Further assume that the array is absolutely regular with exponentially decaying mixing coefficients $\beta_{(n)}(k)\leq \rho^k$ for some $\rho<1$.
If additionally $V_{n,n}:= \var(\sum_{t=1}^n X_{t,n})$ satisfies $\liminf_{n\to \infty} V_{n,n}>0$, then 
\begin{displaymath}
\sum_{t=1}^n \frac{X_{t,n}}{V_{n,n}^{1/2}} \,\stackrel{d}{\longrightarrow}\, Z\sim\mathcal N(0,1).
\end{displaymath}
}
 
\begin{proof}[Proof of Lemma~\ref{p.clt}]
We apply Corollary~1 in \citet{R95}. To this end, we validate conditions (a) and (b) therein. Assumption (a) reads
\begin{displaymath}
\limsup_{n\to\infty} \max_{i\leq n}\frac{\var(\sum_{t=1}^i X_{t,n})}{V_{n,n}} \,<\, \infty.
\end{displaymath}
In view of our assumption $\liminf_{n\to\infty}V_{n,n}>0$, it suffices to show that 
\begin{equation}
\label{eq.rio-a}
\limsup_{n\in\N} \sum_{s,t=1}^n\left|\cov({X_{s,n} }\, , \,{X_{t,n}})\right| \,\leq\, C.
\end{equation}
Recall that the $\beta$-mixing coefficients serve as an upper bound for the $\alpha$-mixing coefficients.
By the covariance inequality for $\alpha$-mixing random variables (see e.g.~\citet[Thm.~3, Sect.~1.2.2]{D94}) we obtain
\begin{displaymath}
\left|\cov({X_{s,n} }\, , \,{X_{t,n}})\right|
\,\leq\,8\, \alpha_{(n)}^X(|s-t|)^{(2+\delta)/\delta}\,\|X_{s,n}\|_{2+\delta}\,\|X_{t,n}\|_{2+\delta}
\end{displaymath}
which then gives \eqref{eq.rio-a} under our moment conditions by the exponential decay of the strong mixing coefficients. 

It remains to check (b) in Corollary~1 in~\citet{R95}, that is 
\begin{equation}
\label{eq.rio-b}
\sum_{t=1}^n \int_0^1 \,\alpha_{(n)}^{-1}\left(\frac{x}{2}\right)\, Q_{t,n}^2(x)\,
\min\left\{1\,,\,\alpha_{(n)}^{-1}\left(\frac{x}{2}\right)\, Q_{t,n}(x)\,\right\}\, dx\ninfty 0. 
\end{equation}
Here, $Q_{t,n}$ denotes the inverse of the survival function $z\mapsto P(|X_{t,n}|\,>\,z\, V_{n,n}^{1/2})$ and
$\alpha_{(n)}^{-1}$ denotes the inverse of $z\mapsto \alpha(\lfloor z\rfloor)$.
First, note that $\alpha_{(n)}^{-1}(x/2)\leq C_1-C_2\,\ln(x)$  for some  $C_1,\, C_2\in (0,\infty)$
in view of the exponential decay of the mixing coefficients. Second, we obtain by Markov's inequality
\begin{displaymath}
P(|X_{t,n}|\,>\,z\, V_{n,n}^{1/2}) \,\leq\, \frac{E|X_{t,n}|^{2+\delta}}{ z^{2+\delta}\, V_{n,n}^{(2+\delta)/2}}
\end{displaymath}
which leads to 
\begin{displaymath}
Q_{t,n}(x)\leq x^{-1/(2+\delta)} \left(\frac{\|X_{t,n}\|_{2+\delta}}{V_{n,n}^{1/2}}\right).
\end{displaymath}
From this, and using the rough estimate
$\min\left\{1,\alpha_{(n)}^{-1}\left(\frac{x}{2}\right)Q_{t,n}(x)\right\}
\leq \alpha_{(n)}^{-\delta/2}\left(\frac{x}{2}\right)Q_{t,n}^{\delta/2}(x)$,
we can deduce that the r.h.s.~of~\eqref{eq.rio-b} can be bounded from above by
\begin{displaymath}
\widetilde C\, \int_0^1 \,|\ln(x)|^{1+\delta/2}\,x^{-(2+\delta/2)/(2+\delta)} \, dx
\,\cdot\,  \sum_{t=1}^n \left(\frac{\|X_{t,n}\|_{2+\delta}}{V_{n,n}^{1/2}}\right)^{2+\delta/2}
\end{displaymath}
which tends to zero as
$\sup_{t\leq n}(\|X_{t,n}\|_{2+\delta})=O(n^{-1/2})$. \hfill\qedhere
\end{proof}

Before we turn to the proof of Proposition~\ref{P3.1} we derive a few useful approximations.
Since $\sigma_0\geq 1$ and
$\ln(\sigma_t)\geq a\ln(\sigma_{t-1}) \,+\, c\, \ln(t)$, we see inductively that under the conditions of Proposition~\ref{P3.1}
\begin{equation}
\label{eq5.1}
\sigma_t \,\geq\, t^c \qquad \forall t\geq 1.
\end{equation} 
In what follows we shall get rid of the cumbersome terms $\ln(X_t+1)$.
Note that the vector $(X_0,\ldots,X_n)$ has the same distribution as
$(\lfloor \sigma_0 Y_0 \rfloor,\ldots,\lfloor \sigma_n Y_n \rfloor)$, where 
$Y_t$ has the same distribution as $Y$ and is independent of $\sigma_0,\ldots,\sigma_t,Y_0,\ldots,Y_{t-1}$. 
To simplify the representation we suppose that $X_t=\lfloor \sigma_t Y_t\rfloor$.

Since
\begin{displaymath}
\ln(\sigma_t Y_t) \leq \ln(X_t+1) \leq \ln(\sigma_t Y_t+1) \leq \ln(\sigma_t(Y_t+1/t^c)),
\end{displaymath}
we obtain
\begin{displaymath}
\ln\big( X_t+1 \big) \,-\, \ln\big( \sigma_t Y_t \big)
\,\leq\, \ln\big( \sigma_t (Y_t+1/t^c) \big) \,-\, \ln\big( \sigma_t Y_t \big)
\end{displaymath}
and hence, for $\rho\leq 2+\delta$,
\begin{eqnarray}
\label{eq5.1}
E\big| \ln\big( X_t+1 \big) \,-\, \ln\big( \sigma_t Y_t \big) \big|^\rho
& \leq & E\big| \ln\big( \sigma_t (Y_t+1/t^c) \big) \,-\, \ln\big( \sigma_t Y_t \big) \big|^\rho \nonumber \\
& = & E\big| \ln\big( (Y_t+1/t^c) \big) \,-\, \ln\big( Y_t \big) \big|^\rho 
\tinfty 0. \qquad
\end{eqnarray}

The above representation allows us to represent $\ln(\sigma_t)$ in a suitable form:
\begin{eqnarray}
\label{eq5.2}
\ln(\sigma_t) 
& = & (a+b)\, \ln(\sigma_{t-1}) \,+\, b\, \ln(Y_{t-1}) \,+\, b\, \big[ \ln(X_{t-1}+1) - \ln(\sigma_{t-1}Y_{t-1}) \big] \,+\, c\, \ln(t) 
\nonumber \\
& = & (a+b)\, \Big\{ (a+b)\, \ln(\sigma_{t-2}) \,+\, b\, \ln(Y_{t-2}) \,+\, b\, \big[ \ln(X_{t-2}+1) - \ln(\sigma_{t-2}Y_{t-2}) \big]
\,+\, c\, \ln(t-1) \Big\} \nonumber \\
& & \quad {} \,+\, b\, \ln(Y_{t-1}) \,+\, b\, \big[ \ln(X_{t-1}+1) - \ln(\sigma_{t-1}Y_{t-1}) \big] \,+\, c\, \ln(t) \nonumber \\
& = & \ldots \,=\, (a+b)^t\, \ln(\sigma_0) \nonumber \\
& & \qquad \,+\, b\, \sum_{s=1}^t (a+b)^{s-1}\, \ln(Y_{t-s}) \nonumber \\
& & \qquad \,+\, b\, \sum_{s=1}^t (a+b)^{s-1}\, \big[ \ln(X_{t-s}+1) \,-\, \ln(\sigma_{t-s}Y_{t-s}) \big] \nonumber \\
& & \qquad \,+\, c\, \Big\{ \sum_{s=0}^{t-1} (a+b)^s\, \ln(t-s) \,-\, \frac{1}{1-a-b}\, \ln(t) \Big\} \nonumber \\
& & \qquad \,+\, \frac{c}{1-a-b}\, \ln(t).
\end{eqnarray}
The first four terms on the right-hand side are stochastically bounded, and the fifth one is dominating.
To see boundedness of the fourth term, note that 
\begin{eqnarray*}
\lefteqn{ \Big| \sum_{s=0}^{t-1} (a+b)^s \ln(t-s) \,-\, \frac{1}{1-a-b} \ln(t) \Big| } \\
& \leq & \Big| \sum_{s=0}^{t-1} (a+b)^s \big[ \ln(t-s) \,-\, \ln(t) \big] \Big|
\,+\, \Big| \frac{1-(a+b)^t}{1-a-b} \ln(t) \,-\, \frac{1}{1-a-b} \ln(t) \Big| \\
& \leq & \sum_{s=0}^{t-1} (a+b)^s \, s \quad + \quad \frac{(a+b)^t}{1-a-b} \, \ln(t).
\end{eqnarray*}

Now we are in a position to approximate the covariance structure of $\ln(X_t+1)$.
Since it follows from \eqref{eq5.1} that $\big\| \ln(X_{t-s}+1) - \ln(\sigma_{t-s}Y_{t-s}) \big\|_2\tinfty 0$,  
we obtain  by Minkowski's inequality
\begin{eqnarray*}
\lefteqn{ \Big\| b \, \sum_{s=1}^t (a+b)^{s-1} \, \big[ \ln(X_{t-s}+1) - \ln(\sigma_{t-s}Y_{t-s}) \big] \Big\|_2 } \\
& \leq & b \, \sum_{s=1}^t (a+b)^{s-1} \, \big\| \ln(X_{t-s}+1) - \ln(\sigma_{t-s}Y_{t-s}) \big\|_2
\,\tinfty\, 0.
\end{eqnarray*}
This means that the limits of the covariances arise from the terms in the fourth from last row in \eqref{eq5.2}:
\begin{equation}
\label{eq5.3}
\cov\big( \ln(\sigma_t), \ln(\sigma_{t-u}) \big) \tinfty \var(\ln(Y)) \, b^2 \, \frac{(a+b)^{u}}{1 \,-\, (a+b)^2}.
\end{equation}
This implies
\begin{eqnarray}
\label{eq5.4}
\var\big( \ln(X_t+1) \big) & \tinfty & \lim_{t\to\infty} \var\big( \ln(\sigma_t) + \ln(Y_t) \big) \nonumber \\
& = & \lim_{t\to\infty} \var\big( \ln(\sigma_t) \big) \quad + \quad \var\big( \ln(Y) \big) \nonumber \\
& = & \var\big( \ln(Y) \big) \, \Big\{ \frac{b^2}{1 - (a+b)^2} \,+\, 1 \Big\}
\end{eqnarray}
and, for $u\geq 1$,
\begin{eqnarray}
\label{eq5.5}
\cov\big( \ln(X_t+1), \ln(X_{t-u}+1) \big)
& \tinfty & \lim_{t\to\infty} \cov\big( \ln(\sigma_t), \ln(\sigma_{t-u}) \big)
\,+\, \lim_{t\to\infty} \cov\big( \ln(\sigma_t), \ln(Y_{t-u}) \big) \nonumber \\
& = & \var\big( \ln(Y) \big) \, \Big\{ \frac{b^2\, (a+b)^{u}}{1-(a+b)^2} \,+\, b\, (a+b)^{u-1} \Big\}.
\end{eqnarray}
In the subsequent proof, \eqref{eq5.4} and \eqref{eq5.5} allow us to identify $\lim_{n\to\infty} V_{n,n}$.

\begin{proof}[Proof of Proposition~\ref{P3.1}]
Let $J_{t,n}=\frac{{\ln(t)}}{\sqrt n\,\ln(n)}\, \left(\ln(X_t+1)-\, E\ln(X_t+1) \right)$.  
By Lemma~\ref{LA.3} $$\sum_{t=1}^n \ln^2(t)/(n\, \ln^2(n))\ninfty 1.$$
Therefore, it remains to verify
\begin{displaymath}
\sum_{t=1}^n J_{t,n} \,\stackrel{d}{\longrightarrow}\, Z\sim \mathcal N(0,\sigma^2).
\end{displaymath}
To this end, we apply Lemma~\ref{p.clt} with $X_{t,n}:=J_{t,n}$.  Note that we obtain
\begin{eqnarray*}
\sup_{t\leq n }\|J_{t,n}\|_{2+\delta}
& \leq & O(n^{-1/2}) \, \sup_{t\leq n,\, n\in\N}  \|\ln(X_t+1)-\, E\ln(X_t+1)\|_{2+\delta} \\
& \leq & O(n^{-1/2}) \, \sup_{t\leq n,\, n\in\N} \Big[
	\|\ln(X_t+1)-\,  \ln(\sigma_t\, Y_t)\|_{2+\delta}
\\
&	& \qquad {}	+ \|\ln(\sigma_t\, Y_t)-\, E\ln(\sigma_t\, Y_t )\|_{2+\delta}+ \|E[\ln(\sigma_t\, Y_t)-\, \ln(X_t+1)]\|_{2+\delta}
	\Big].
\end{eqnarray*}
Uniform boundedness of the first and the last summand can be deduced from~\eqref{eq5.1}.
Uniform boundedness of the middle term follows from Remark~\ref{R3} together with \eqref{eq5.2}
since the latter gives
\begin{eqnarray}
\label{eq.fin-mom}
\lefteqn{ \|\ln(\sigma_t\, Y_t)-\, E\ln(\sigma_t\, Y_t )\|_{2+\delta} } \nonumber \\
& \leq & \|\ln(\sigma_t)-\, E\ln(\sigma_t)\|_{2+\delta} \,+\, \|\ln(Y_t)-\, E\ln(Y_t )\|_{2+\delta} \nonumber \\
& \leq & (a+b)^t \, \| \ln(\sigma_0)-\, E\ln(\sigma_0)\|_{2+\delta} \nonumber \\
& & {} \,+\, \sum_{s=1}^t (a+b)^{s-1}\, b\, \big\| [ \ln(X_{t-s}+1) \,-\, \ln(\sigma_{t-s}) ]-E[ \ln(X_{t-s}+1)
\,-\, \ln(\sigma_{t-s}) ]\big\|_{2+\delta} \nonumber \\
& & {} \,+\, \| \ln(Y) \,-\, E\ln(Y) \|_{2+\delta}.
\end{eqnarray}
It remains to derive $\liminf_{n\to\infty } V_{n,n}$.
It follows from \eqref{eq5.1} and \eqref{eq.fin-mom} that the moments of order $2+\delta$ of the random variables $\ln(X_t+1)$
are bounded. Hence, we obtain from the mixing property that
$\sup_{t}|\cov(\ln(X_{t}+1), \ln(X_{t+h}+1))|\leq C\, \rho^{|h|}$ for some $\rho\in (0,1)$.
Furthermore, Lebesgue's theorem, Cauchy's limit theorem, Lemma~\ref{LA.3}, \eqref{eq5.4}  and \eqref{eq5.5} give
\begin{eqnarray*}
\lim_{n\to\infty} V_{n,n}
& = & \sum_{h\in\Z}\lim_{n\to\infty}\frac{1}{n}\sum_{t=1}^n\,\frac{\ln(t+h)\,\ln(t)}{\ln^2(n)}\,\cov(\ln(X_{t}+1), \ln(X_{t+h}+1))
\,\1_{1\leq t+h\leq n} \\
& = & \var\big( \ln(Y) \big) \,\Big\{1\,+\,\sum_{h\in\Z}   \frac{b^2\, (a+b)^{|h|}}{1-(a+b)^2} \,+\, 2\sum_{h\in\N} b\, (a+b)^{h-1}\, \Big\}   \\
\lefteqn{ \;=\, \var\big( \ln(Y) \big) \,  \frac{(1-a)^2}{(1-(a+b))^2} \,>\, 0. \hspace*{6.13cm}\qedhere}
\end{eqnarray*} 
\end{proof}
%\pagebreak

%%%%%%%%%%%%%%%%%%%%%%%%%%%%%%%%%%%%%%%%%%%%%%%%%%%%%%%%%%%%%
\section{Proof of Proposition~\ref{P3.2}}
\label{S6neu}
%%%%%%%%%%%%%%%%%%%%%%%%%%%%%%%%%%%%%%%%%%%%%%%%%%%%%%%%%%%%%

	\begin{proof}[Proof of Proposition~\ref{P3.2}]{}\
		\begin{itemize}
			\item[(i)\quad]
			It follows from Proposition~\ref{P3.1} that
			\begin{equation}
				\label{pp32.1}
				\sup_{x\in\R} \big| P(T_n\leq x) \,-\, \Phi(x/\sigma) \big| \,\ninfty\, 0.
			\end{equation}
			Since $\sigma_n^2:=\var(T_n)\ninfty\sigma^2$ and since $T_n^*$ conditioned on $X_1,\ldots,X_n$
			is by construction normally distributed, it only remains to show, 
			for $\sigma_n^{*2}:=\var(T_n^*\mid X_1,\ldots,X_n)$, that
			\begin{equation}
				\label{pp32.2}
				\big| \sigma_n^{*2} \,-\, \sigma_n^2 \big| \,\stackrel{P}{\longrightarrow}\, 0.
			\end{equation}
			Then (i)  follows from (\ref{pp32.1}) and (\ref{pp32.2}).\\
			Let $Z_t=w_t[\ln(X_t+1)-E\ln(X_t+1)]$ and $Z_t'=w_t[E\ln(X_t+1)-\widehat{m}_n(t)]$.
			Then $T_n=\sum_{t=1}^n Z_t$ and $T_n^*=\sum_{t=1}^n (Z_t+Z_t') W_t^*$.
			Note that it follows from (\ref{eq3.3}), (\ref{eq3.4}) and $w_t=O(n^{-1/2})$ that
			$\sum_{t=1}^n E\big[(Z_t')^2\big]=O\big(n^{-1}\,(n/N_n + N_n(\ln(n))^2)\big)$.
			Since $\sum_{r=-(n-1)}^{n-1}|\rho(r/l_n)|=O(l_n)$, we obtain that
		  \begin{eqnarray*}
				E\Big[ \Big( \sum_{t=1}^n Z_t' W_t^* \Big)^2 \Big]
				& = & E\Big[ E\Big( \Big( \sum_{t=1}^n Z_t' W_t^* \Big)^2 \mid X_1,\ldots,X_n \Big) \Big] \\
				& = & E\Big[ \sum_{s,t=1}^n Z_s' Z_t' \rho(|s-t|/l_n) \Big] \\
	      & = & \sum_{r=-(n-1)}^{n-1} \big| \rho(r/l_n) \big| \, \sum_{\max\{1,1-r\}\leq s\leq \min\{n,n-r\}} \big| E\big[ Z_s' Z_{s+r}' \big] \big| \\
				& = & O\big( \frac{l_n}{n} \, \big( \frac{n}{N_n} \,+\, N_n (\ln(n))^2 \big) \big) \,=\, o(1).
			\end{eqnarray*}
			Next we compare $\var(T_n)$ and $\var(\sum_t Z_t W_t^*\mid X_1,\ldots,X_n)$.
			We obtain
			\begin{eqnarray}
				\label{pp32.3}
				\lefteqn{ E\big[T_n^2\big] \,-\, E\big((\sum_t Z_t W_t^*)^2 \mid X_1,\ldots,X_n\big) } \nonumber \\
				& = & \sum_{s,t} E\big[ Z_s Z_t \big] \,-\, \sum_{s,t} Z_s Z_t \rho(|s-t|/l_n) \nonumber \\
				& = & \sum_{s,t} E\big[ Z_s Z_t (1 \,-\, \rho(|s-t|/l_n)) \big] \nonumber \\
				& & {} \,-\, \sum_{s,t} \big( Z_s Z_t \,-\, E[Z_s Z_t] \big) \rho(|s-t|/l_n).
			\end{eqnarray}
			The last but one term converges to zero by majorized convergence.
			We conclude from (\ref{eq5.1}) and (\ref{eq5.2}) that
			\begin{displaymath}
				\sup_t \big\{ E| \ln(X_t+1) - E\ln(X_t+1)|^{4+\delta} \big\} \,<\, \infty.
			\end{displaymath}
			This and the exponential decay of the mixing coefficients $\beta^X(r)$ imply,
			for $1\leq s\leq t\leq u\leq v\leq n$ and $r=\max\{t-s,u-t,v-u\}$, that
			\begin{displaymath}
				\big| \mbox{cum}(Z_s,Z_t,Z_u,Z_v) \big| \,=\, O\big( \rho^r \, n^{-2} \big)
			\end{displaymath}
			for some $\rho<1$, where\\
			%\begin{displaymath}
			$\mbox{cum}(Z_s,Z_t,Z_u,Z_v)=E[Z_sZ_tZ_uZ_v]\!-\!E[Z_sZ_t]E[Z_uZ_v]\!-\!E[Z_sZ_u]E[Z_tZ_v]\!-\!E[Z_sZ_v]E[Z_tZ_u]$
			%\end{displaymath}
			denotes the joint cumulant of $Z_s$, $Z_t$, $Z_u$ and $Z_v$.
			Hence,
			\begin{displaymath}
				\sum_{s,t,u,v=1}^n \big| \mbox{cum}(Z_s,Z_t,Z_u,Z_v) \big|
				\,\leq\, 4!\, \sum_{1\leq s\leq t\leq u\leq v\leq n} \big| \mbox{cum}(Z_s,Z_t,Z_u,Z_v) \big|
				\,=\, O(1/n).
			\end{displaymath}
			Furthermore, it follows from (\ref{eq5.4}) and (\ref{eq5.5}) that $\cov(Z_s,Z_t)=O(n^{-1}\, (a+b)^{|s-t|})$.
			Therefore we obtain for the second term on the right-hand side of (\ref{pp32.3}) that
			\begin{eqnarray*}
				\lefteqn{ E\Big[ \big( \sum_{s,t} \big( Z_s Z_t \,-\, E[Z_s Z_t] \big) \rho(|s-t|/l_n) \big)^2 \Big] } \\
				& = & \sum_{s,t,u,v=1}^n \rho(|s-t|/l_n) \, \rho(|u-v|/l_n) \, \mbox{cum}(Z_s,Z_t,Z_u,Z_v) \\
				& & {} \,+\, 2\, \sum_{s,t,u,v=1}^n \rho(|s-t|/l_n) \, \cov(Z_t,Z_u) \, \rho(|u-v|/l_n) \, \cov(Z_v,Z_s) \\
				& = & O\Big( \frac{1}{n} \Big) \,+\, O\Big( \frac{l_n}{n} \Big).
			\end{eqnarray*}
			\item[(ii)\quad]
			Follows from (i). \hfill\qedhere
		\end{itemize}
	\end{proof}

	{\lem
		\label{L.monotone}
		Suppose that (\ref{1.1a}) and (\ref{eq3.0}) are fulfilled and that $\sigma_0=1$.
		Then the function $t\mapsto E\ln(X_t+1)$ is monotonously increasing.
	}

	\begin{proof}[Proof of Lemma~\ref{L.monotone}]
		We construct two versions $(\widetilde{\sigma}_t)_{t\in\N_0}$ and $(\widetilde{\sigma}_t')_{t\in\N_0}$
		of the intensity process such that $\widetilde{\sigma}_0=\widetilde{\sigma}_0'=1$ and
		\begin{equation}
			\label{lm.1}
			\widetilde{\sigma}_t \,\leq\, \widetilde{\sigma}_{t+1}' \qquad \forall t\in\N
		\end{equation}
		holds with probability~1.
		This will be achieved by feeding the first process with innovations $Y_1,Y_2,\ldots$
		while the second one is fed with $Y_0,Y_1,Y_2,\ldots$, where $(Y_t)_{t\in\N_0}$ is a sequence
		of independent random variables such that $Y_t\stackrel{d}{=}Y$, i.e.~we define for $t\geq 1$
		\begin{eqnarray*}
			\ln(\widetilde{\sigma}_t) & = & a\, \ln( \widetilde{\sigma}_{t-1} )
			\,+\, b\, \ln( \lfloor \widetilde{\sigma}_{t-1} Y_t \rfloor \,+\, 1 ) \,+\, c\, \ln(t) \\
			\mbox{and} & & \\
			\ln(\widetilde{\sigma}_t') & = & a\, \ln( \widetilde{\sigma}_{t-1}' )
			\,+\, b\, \ln( \lfloor \widetilde{\sigma}_{t-1}' Y_{t-1} \rfloor \,+\, 1 ) \,+\, c\, \ln(t).
		\end{eqnarray*} 
		Then
		\begin{eqnarray*}
			\ln( \widetilde{\sigma}_1' ) & = & a\, \ln( \widetilde{\sigma}_0' )
			\,+\, b\, \ln( \lfloor \widetilde{\sigma}_0' Y_0 \rfloor \,+\, 1 ) \,+\, c\, \ln(1) \\
			& \geq & 0 \,=\, \ln( \widetilde{\sigma}_0 ).
		\end{eqnarray*}
		Suppose now that, for $t\geq 1$, $\ln(\widetilde{\sigma}_t')\geq \ln(\widetilde{\sigma}_{t-1})$. Then
		\begin{eqnarray*}
			\ln( \widetilde{\sigma}_{t+1}' ) & = & a\, \ln( \widetilde{\sigma}_t' )
			\,+\, b\, \ln( \lfloor \widetilde{\sigma}_t' Y_t \rfloor \,+\, 1 ) \,+\, c\, \ln(t+1) \\
			& > & a\, \ln( \widetilde{\sigma}_{t-1} )
			\,+\, b\, \ln( \lfloor \widetilde{\sigma}_{t-1} Y_t \rfloor \,+\, 1 ) \,+\, c\, \ln(t) \\
			& = & \ln( \widetilde{\sigma}_t ),
		\end{eqnarray*}
		which proves that (\ref{lm.1}) holds for all $t\in\N$ and with probability~1.
		This implies that
		\begin{displaymath}
		\hspace*{1cm}
			E\ln( X_t+1 ) \,=\, E\ln( \lfloor \widetilde{\sigma}_t Y_t \rfloor + 1 )
			\,<\, E\ln( \lfloor \widetilde{\sigma}_{t+1}' Y_t \rfloor + 1 ) \,=\, E\ln( X_{t+1}+1 ). \hspace*{1cm}\qedhere
		\end{displaymath}
	\end{proof}

%%%%%%%%%%%%%%%%%%%%%%%%%%%%%%%%%%%%%%%%%%%%%%%%%%%%%%%%%%%%%%%%%%%%%%%%%%%%%%%
%%%%%%%%%%%%%%%%%%%%%%%%%%%%%%%%%%%%%%%%%%%%%%%%%%%%%%%%%%%%%%%%%%%%%%%%%%%%%%%
\section{A few auxiliary results}
\label{S6}
%%%%%%%%%%%%%%%%%%%%%%%%%%%%%%%%%%%%%%%%%%%%%%%%%%%%%%%%%%%%%%%%%%%%%%%%%%%%%%%

In this section we collect a few technical results which contribute to the proof of our main result.

{\lem
\label{L.tv}
Suppose that $Y$ is a non-negative random variable with a density~$p$.
\begin{itemize}
\item[(i)] If $p$ is monotonously non-increasing on $[0,\infty)$, then
\begin{displaymath}
d_{TV}\big(P^{\lfloor \sigma Y\rfloor}, P^{\lfloor \sigma' Y \rfloor}\big) \,\leq\, \big| \ln(\sigma) \,-\, \ln(\sigma') \big|
\qquad \forall \sigma,\sigma'>0.
\end{displaymath}
\item[(ii)] If $p$ is everywhere differentiable on $(0,\infty)$ and $\int_0^\infty x|p'(x)| \, dx<\infty$, then
\begin{displaymath}
d_{TV}\big(P^{\lfloor \sigma Y\rfloor}, P^{\lfloor \sigma' Y \rfloor}\big)
\,\leq\, \big| \ln(\sigma) \,-\, \ln(\sigma') \big|\; \Big\{ 1 \,+\, \int_0^\infty x|p'(x)| \, dx \Big\}/2.
\qquad \forall \sigma,\sigma'>0.
\end{displaymath}
\end{itemize}
}

\begin{proof}
Note that the total variation distance 
$d_{TV}\big( P_\sigma, P_{\sigma'} \big)=(1/2)\sum_{k=0}^\infty|P_\sigma(\{k\})-P_{\sigma'}(\{k\})|$
between~$P_\sigma$ and~$P_{\sigma'}$ can be bounded from above by
\begin{eqnarray}
\label{1.2}
d_{TV}\big( P_\sigma, P_{\sigma'} \big)
& = & 1 \,-\, \sum_{k=0}^\infty P_\sigma(\{k\}) \wedge P_{\sigma'}(\{k\}) \nonumber \\
& \leq & 1 \,-\, \int_0^\infty \frac{1}{\sigma} p(x/\sigma) \wedge \frac{1}{\sigma'} p(x/\sigma') \, dx \nonumber \\
& = & \frac{1}{2} \int_0^\infty \big| \frac{1}{\sigma} p\big( \frac{x}{\sigma} \big)
\,-\, \frac{1}{\sigma'} p\big( \frac{x}{\sigma'} \big) \big|\, dx, 
\end{eqnarray}
which is just the total variation distance between the distributions of $\sigma Y$ and $\sigma' Y$.

Let $\sigma<\sigma'$.
If $p$ is non-increasing on $[0,\infty)$, then
\begin{eqnarray*}
1 \,-\, \int_0^\infty \frac{1}{\sigma} p(x/\sigma) \wedge \frac{1}{\sigma'} p(x/\sigma') \, dx
& \leq & 1 \,-\, \int_0^\infty \frac{1}{\sigma'} \, p\big( \frac{x}{\sigma} \big) \, dx \\
& = & 1 \,-\, \frac{\sigma}{\sigma'} \,\leq\, \int_\sigma^{\sigma'} \frac{1}{x} \, dx
\,=\, \ln(\sigma') \,-\, \ln(\sigma).
\end{eqnarray*}
If $p$ is everywhere differentiable on $(0,\infty)$ and $\int_0^\infty x|p'(x)|\, dx<\infty$, then
\begin{eqnarray*}
\frac{1}{2} \int_0^\infty \big| \frac{1}{\sigma} p\big( \frac{x}{\sigma} \big)
\,-\, \frac{1}{\sigma'} p\big( \frac{x}{\sigma'} \big) \big| \, dx
& \leq & \frac{1}{2} \int_0^\infty \int_\sigma^{\sigma'} \Big|\frac{\partial}{\partial u} \Big( \frac{1}{u}
p\big( \frac{x}{u} \big) \Big) \Big| \, du \, dx \\
& \leq & \frac{1}{2} \int_0^\infty \int_\sigma^{\sigma'} \frac{1}{u^2} p\big( \frac{x}{u} \big)
\,+\, \frac{x}{u^3} \big| p'\big( \frac{x}{u} \big) \big| \, du \, dx \\
& = & \frac{1}{2} \, \int_\sigma^{\sigma'} \frac{1}{u} \bigg\{ \int_0^\infty \frac{1}{u} p\big( \frac{x}{u} \big)
\,+\, \frac{x}{u^2} \big| p'\big( \frac{x}{u} \big) \big| \, dx \bigg\} \, du \\
& = & \frac{ \ln(\sigma') \,-\, \ln(\sigma) }{ 2 } 
\Big\{ 1 \,+\, \int_0^\infty x\big| p'(x) \big| \, dx \Big\},
\end{eqnarray*}
which completes the proof. \hfill\qedhere
\end{proof}

{\lem
\label{LA.1}
Suppose that $Y$ is a non-negative random variable with a bounded density~$p$ and  $E[\ln^+(Y)]<\infty$.
Then
\begin{displaymath}
E\big| \ln( \lfloor \sigma Y \rfloor + 1 ) \,-\, \ln( \sigma+1 ) \big| \,\leq\, \|p\|_\infty  \,+\, E[\ln^+(Y)]
\qquad \forall \sigma>0.
\end{displaymath}
}

\begin{proof}
Using
\begin{eqnarray*}
\sum_{k=0}^{\lfloor \sigma \rfloor} \big( \ln(\sigma+1) \,-\, \ln(k+1) \big)
& = & \big( \lfloor \sigma \rfloor + 1 \big) \, \ln(\sigma+1) \,-\, \sum_{k=1}^{\lfloor \sigma \rfloor} \ln(k+1) \\
& \leq & \big( \lfloor \sigma \rfloor + 1 \big) \, \ln(\sigma+1) \,-\, \int_1^{\lfloor \sigma \rfloor+1} \ln(y) \, dy \\
& = & \big( \lfloor \sigma \rfloor + 1 \big) \, \ln(\sigma+1)
\,-\, \big[ x\ln(x) \,-\, x \big]_1^{\lfloor \sigma \rfloor +1} \\
& = & \big( \lfloor \sigma \rfloor + 1 \big) \, \underbrace{\big( \ln(\sigma+1) \,-\, \ln(\lfloor \sigma \rfloor + 1) \big)}
_{\leq\, (\sigma-\lfloor\sigma\rfloor)/(\lfloor\sigma \rfloor+1)}
\,+\, \lfloor \sigma \rfloor \\
& \leq & \sigma,
\end{eqnarray*}
we obtain that
\begin{eqnarray*}
\lefteqn{ E\big[ |\ln(\lfloor\sigma Y\rfloor+1)\,-\,\ln(\sigma+1)|\;\1(\lfloor\sigma Y\rfloor\leq\sigma) \big] } \\
& = & \sum_{k=0}^{\lfloor \sigma \rfloor} \big( \ln(\sigma+1) \,-\, \ln(k+1) \big)
\, P\big(\sigma Y\in [k,k+1) \big) \,\leq\, \|p\|_\infty.
\end{eqnarray*}
Furthermore, we have that
\begin{eqnarray*}
\lefteqn{ E\big[ |\ln(\lfloor\sigma Y\rfloor+1)\,-\,\ln(\sigma+1)|\;\1(\lfloor\sigma Y\rfloor>\sigma) \big] } \\
& \leq & E\big[ (\ln(\sigma Y +1)\,-\,\ln(\sigma +1))\;\1(Y\geq 1) \big] \\
& \leq & E\big[ (\ln(\sigma Y)\,-\,\ln(\sigma))\;\1(Y\geq 1) \big] \\
& = & E\big[ \ln^+(Y) ],
\end{eqnarray*}
which completes the proof. \hfill\qedhere
\end{proof}

{\lem
\label{LA.2}
Suppose that $Y$ is a non-negative random variable with a continuous probability density~$p$ satisfying
$\gamma:=\int_0^\infty \sup\{p(y)\colon\, y\geq x\}\, dx<\infty$.
Then the function \mbox{$\sigma\mapsto E\ln(\lfloor \sigma Y\rfloor+1)$} is differentiable and
\begin{displaymath}
\frac{d}{d\sigma} E\ln\big( \lfloor \sigma Y \rfloor \,+\, 1 \big) \,\leq\, \frac{\gamma}{\sigma} \qquad \forall \sigma>0.
\end{displaymath}
(If $p$ is monotonously non-increasing, then $\gamma=1$.)
}

\begin{proof}
First note that
\begin{displaymath}
E \ln\big( \lfloor \sigma Y \rfloor \,+\, 1 \big)
\,=\, \sum_{k=1}^\infty \big( \ln(k+1) \,-\, \ln(k) \big) P\big( \lfloor \sigma Y \rfloor \geq k \big)
\,=\, \sum_{k=1}^\infty \big( \ln(k+1) \,-\, \ln(k) \big) P\big( \sigma Y \geq k \big).
\end{displaymath}
To prove differentiability we consider corresponding difference quotients. Let $g(x):=\sup\{p(y)\colon\, y\geq x\}$.
For \mbox{$0<\epsilon<\sigma$}, we have that
\begin{eqnarray*}
\frac{ P\big( (\sigma+\epsilon)Y \geq k \big) \,-\, P\big( \sigma Y \geq k \big) }{ \epsilon }
& = & \frac{1}{\epsilon} \, \int_{k/(\sigma+\epsilon)}^{k/\sigma} p(y) \, dy \\
& \leq & \frac{1}{\epsilon} \, k \big( \frac{1}{\sigma} \,-\, \frac{1}{\sigma+\epsilon} \big) \, g\big( \frac{k}{\sigma+\epsilon} \big) \\
& \leq & \frac{k}{\sigma^2} \, g\big( \frac{k}{2\sigma} \big) 
\,\leq\, \frac{k}{\sigma} \, \int_{k-1}^k \frac{1}{\sigma} g\big( \frac{y}{2\sigma} \big) \, dy \,=:\, h_1(k),
\end{eqnarray*}
and, for $-\sigma/2<\epsilon<0$, we obtain that 
\begin{eqnarray*}
\frac{ P\big( (\sigma+\epsilon)Y \geq k \big) \,-\, P\big( \sigma Y \geq k \big) }{ \epsilon }
& = & -\frac{1}{\epsilon} \, \int_{k/\sigma}^{k/(\sigma+\epsilon)} p(y) \, dy \\
& \leq & \frac{1}{-\epsilon} \, k \big( \frac{1}{\sigma+\epsilon} \,-\, \frac{1}{\sigma} \big) \, g\big( \frac{k}{\sigma} \big) \\
& \leq & \frac{k}{\sigma^2/2} \, g\big( \frac{k}{\sigma} \big)
\,\leq\, \frac{k}{\sigma} \, \int_{k-1}^k \frac{1}{\sigma/2} \, g\big( \frac{y}{\sigma} \big) \, dy \,=:\, h_2(k).
\end{eqnarray*}
In both cases, the corresponding quantities $h_1(k)$ and $h_2(k)$ do not depend on the value of $\epsilon$, and we obtain that 
\begin{displaymath}
\sum_{k=1}^\infty \underbrace{ \big( \ln(k+1) \,-\, \ln(k) \big) }_{\leq\, 1/k} h_i(k)
\,\leq\, \frac{2}{\sigma} \, \int_0^\infty g(y)\, dy \,<\, \infty,
\end{displaymath}
which allows us to invoke Lebesgue's theorem on dominated convergence.
Since $p$ is continuous, we have 
\begin{displaymath}
\lim_{\epsilon\to 0} \frac{ P\big( (\sigma+\epsilon)Y \geq k \big) \,-\, P\big( \sigma Y \geq k \big) }{ \epsilon } 
\,=\, \frac{k}{\sigma^2} \, p\big( \frac{k}{\sigma} \big),
\end{displaymath}
which implies that
\begin{eqnarray*}
\frac{d}{d\sigma} E\ln\big( \lfloor \sigma Y \rfloor \,+\, 1 \big)
& = & \sum_{k=1}^\infty \underbrace{k\, \big( \ln(k+1) \,-\, \ln(k) \big)}_{\leq\, 1} \;
\frac{1}{\sigma^2} p\big( \frac{k}{\sigma} \big) \\
\lefteqn{ \;\leq\, \frac{1}{\sigma} \sum_{k=1}^\infty \int_{k-1}^k \frac{1}{\sigma} g(y/\sigma) \, dx
\,\leq\, \frac{\gamma}{\sigma}. \hspace*{2.98cm}\qedhere}
\end{eqnarray*}
\end{proof}

{\lem
\label{LA.3}
It holds, for fixed $h\in\Z$,
\begin{displaymath}
\sum_{t\colon\, 1\leq t,t+h\leq n} \ln(t+h)\, \ln(t) \,=\, n\, \big( \ln(n) \big)^2 \,+\, O\big( n\, \ln(n) \big).
\end{displaymath}
}

\begin{proof}
Let, w.l.o.g., $h\geq 0$. Then
\begin{displaymath}
\sum_{t\colon\, 1\leq t,t+h\leq n} \ln(t+h)\, \ln(t) \\
\,=\, (n-h)\, \big( \ln(n) \big)^2 \,-\, \sum_{t=1}^{n-h} \big( \ln(n) \big)^2 \,-\, \ln(t+h)\, \ln(t).
\end{displaymath}
Since $\ln(s+1)-\ln(s)\leq s^{-1}$ for $s\geq 1$    we can estimate the right-hand side by
\begin{eqnarray*}
\lefteqn{ \sum_{t=1}^{n-h} \big( \ln(n) \big)^2 \,-\, \ln(t+h)\, \ln(t) } \\
& \leq & \sum_{t=1}^{n-h} \big( \ln(n) \big)^2 \,-\, \big( \ln(t) \big)^2 \\
& = & \sum_{t=1}^{n-h} \big( \ln(n) \big)^2 \,-\, \big( \ln(n-1) \big)^2 \,+\, \big( \ln(n-1) \big)^2
\,-\, \cdots \,-\, \big( \ln(t+1) \big)^2 \,+\, \big( \ln(t+1) \big)^2 \,-\, \big( \ln(t) \big)^2 \\
& = & \sum_{s=1}^{n-1} \min\{s,n-h\} \, \left[\big( \ln(s+1) \big)^2 \,-\, \big( \ln(s) \big)^2 \right]\\
& \leq & \sum_{s=1}^{n-1} \ln(s+1) \,+\, \ln(s)\\
&  \leq& 2\, \sum_{t=1}^n \ln(t)
\,=\, O\big( n\, \ln(n) \big),
\end{eqnarray*}
which completes the proof. \hfill\qedhere
\end{proof}

%%%%%%%%%%%%%%%%%%%%%%%%%%%%%%%%%%%%%%%%%%%%%%
%% Single Appendix:                         %%
%%%%%%%%%%%%%%%%%%%%%%%%%%%%%%%%%%%%%%%%%%%%%%
%\begin{appendix}
%\section*{???}%% if no title is needed, leave empty \section*{}.
%\end{appendix}
%%%%%%%%%%%%%%%%%%%%%%%%%%%%%%%%%%%%%%%%%%%%%%
%% Multiple Appendixes:                     %%
%%%%%%%%%%%%%%%%%%%%%%%%%%%%%%%%%%%%%%%%%%%%%%
%\begin{appendix}
%\section{???}
%
%\section{???}
%
%\end{appendix}

%%%%%%%%%%%%%%%%%%%%%%%%%%%%%%%%%%%%%%%%%%%%%%
%% Support information, if any,             %%
%% should be provided in the                %%
%% Acknowledgements section.                %%
%%%%%%%%%%%%%%%%%%%%%%%%%%%%%%%%%%%%%%%%%%%%%%
\begin{acks}[Acknowledgments]
This work was partly funded by Project ``EcoDep'' PSI-AAP2020 -- 0000000013 and by "ProBe-Pro-Oberfranken" FP01054.
The authors thank the editors and two anonymous referees for their valuable comments that led to a significant improvement of the paper.
\end{acks}

\bibliographystyle{harvard}

\end{document}